\documentclass[12pt, reqno]{amsart}
\usepackage{amsfonts,amssymb,graphicx,enumerate,youngtab}
\usepackage{fourier}

\numberwithin{equation}{section}

\theoremstyle{definition}

\theoremstyle{plain}
\newtheorem{theorem}{Theorem}[section]

\newtheorem{corollary}[theorem]{Corollary}
\newtheorem{lemma}[theorem]{Lemma}

\newcommand{\CC}{\mathbf{C}}
\newcommand{\HH}{\mathbf{H}}
\newcommand{\RR}{\mathbf{R}}

\newcommand{\ZZ}{\mathbf{Z}}

\begin{document}

\title{The $n!$ theorem via quaternionic rotation of Cherednik algebras}

\author{Stephen Griffeth}

\address{Instituto de Matem\'aticas \\ Universidad de Talca \\ sgriffeth@utalca.cl}

\maketitle

\begin{abstract}
We give a proof of the $n!$ theorem based on the representation theory of Cherednik algebras and an analog of the rotation of complex structures arising from an action of the quaternions.
\end{abstract}

\section{Introduction}

\subsection{The $n!$ theorem} Let $n$ be a positive integer, let $S=\CC[x_1,\dots,x_n,y_1,\dots,y_n]$ be the ring of polynomials in $2n$ variables, and let $\lambda=(\lambda_1 \geq \lambda_2 \geq \cdots \geq \lambda_\ell>0)$ be an integer partition of $n$, identified with its \emph{diagram}
$$\lambda=\{(a,b) \in \ZZ_{\geq 0}^2 \ | \ 0 \leq b \leq \ell-1, \  0 \leq a \leq \lambda_{b+1}-1\},$$ which is a set of size $n$. Enumerating so that
$$\lambda=\{(a_1,b_1),(a_2,b_2),\dots,(a_n,b_n) \}$$ we define the polynomial $\Delta_\lambda \in S$ by
$$\Delta_\lambda=\mathrm{det}(x_i^{a_j} y_i^{b_j}),$$ and consider the map from $S$ to $S$ defined by sending the polynomial $p \in S$ to 
$$p(\partial) \cdot \Delta_\lambda,$$ where $p(\partial)$ is the differential operator obtained from $p$ by replacing $x_i$ and $y_i$ by the corresponding partial derivatives $\partial_{x_i}$ and $\partial_{y_i}$. The \emph{Garsia-Haiman module} $\mathrm{GH}_\lambda$ is the quotient of $S$ by the kernel of this map. The $n!$ theorem is
\begin{theorem}[Haiman]
The dimension of $\mathrm{GH}_\lambda$ is $n!$. 
\end{theorem} Haiman (in \cite{Hai1}) deduced Macdonald's positivity conjecture from this theorem, as well as an isomorphism of algebraic varieties between the Hilbert scheme of $n$ points in the plane $\mathrm{Hilb}_n(\CC^2)$ and the $S_n$-Hilbert scheme $S_n\mathrm{-Hilb}(\CC^{2n})$ of $\CC^{2n}$, and then by using it proved (in \cite{Hai2}) that the dimension of the diagonal coinvariant ring (the quotient of $S$ by the ideal generated by $S_n$-invariants of positive degree) is $(n+1)^{n-1}$. 

\subsection{Zero fiber rings and quaternions} 

In addition to his work on the symmetric group, Haiman conjectured that for an irreducible real reflection group of rank $n$ the dimension of the diagonal coinvariant ring is at least $(h+1)^n$, where $h$ is the Coxeter number. Gordon \cite{Gor} proved this, using the representation theory of Cherednik algebras. But even in the case of the symmetric group, Gordon's result only implies the (correct) lower bound for the dimension, and the underlying reason that the technique works remains mysterious. This mystery and the so-far-unmet expectation that Cherednik algebras might be able to contribute more than just lower bounds to this circle of ideas seem to me to be some of the most important open problems in the field.

In recent years, we have extended the connection between diagonal coinvariant rings and Cherednik algebras, first to complex reflection groups (in \cite{Gri2}) and more recently (and partly conjecturally) to quaternionic reflection groups (in \cite{CaGr}; the paper \cite{GrGu} observes some striking parallels between the numerology of quaternionic reflection groups and that of complex reflection groups). The extension of Haiman's conjecture on diagonal coinvariant rings of real reflection groups to the zero fiber rings of quaternionic reflection groups suggests that quaternionic structures should be involved in some way in the phenomena first observed by Haiman, even in the case of the symmetric group $S_n$ acting on $(\CC^2)^{\oplus n} \cong \HH^n$.

Taking this philosophy seriously leads, in the case of the symmetric group, to the ideas and methods in this paper. Somewhat more specifically, we use a rotation of the Cherednik algebra action on an irreducible module by a quaternionic action together with the observation of \cite{Gri1} showing that the Garsia-Haiman modules are quotients of associated gradeds of irreducible representations of the Cherednik algebra to give a new proof of the $n!$ theorem. We hope and believe that the proof we present here of the $n!$ theorem is part of a beautiful landscape of geometric (both algebraic and metric), algebraic, and combinatorial phenomena around quaternionic reflection groups, mostly still unexplored. 

\subsection{Note on the exposition} 

This paper is written for experts, and is therefore quite terse in places. In particular, we omit algebraic manipulations we believe are routine. We are preparing an expanded version, based on lectures the author gave in Montreal in 2017 together with an ongoing seminar in Talca on the $n!$ theorem. The goal of the expanded notes is to provide a self-contained proof of the $n!$ theorem, explaining (for a target audience of graduate students in mathematics) all the results on representation theory of Cherednik algebras, the Kempf-Ness theorem, and the algebraic combinatorics necessary for it.

\subsection{Use of AI}

This is the first paper I have written in which AI played a significant role in the development of the ideas. Specifically, the strategy here and all the writing is mine, but many of the tactics and preliminary computations were developed in a lengthy sequence of interactions with Claude Opus 5.5 and one interaction (during a promotional period in which it was freely available) with Claude Fable 5.2. Had similar interactions occurred with a human then co-authorship would certainly have been appropriate. I very much doubt that this paper would exist in its current form without those interactions.

\section{Minimal metrics and rotation of the $H$-action}

\subsection{Metrics} Let $V$ be a finite-dimensional $\CC$-vector space. A \emph{Hermitian form} on $V$ is a function $h:V \times V \to \CC$ such that
\begin{itemize}
\item[(a)] $h(u,a_1 v_1+a_2 v_2)=a_1 h(u, v_1)+a_2h(u,v_2)$ for all $u,v_1,v_2 \in V$ and $a_1,a_2 \in \CC$, and
\item[(b)] $h(u,v)=\overline{h(v,u)}$ for all $u,v \in V$.
\end{itemize} A Hermitian form $h$ is \emph{positive definite} if 
$$h(u,u) > 0 \quad \hbox{for all $u \in V \setminus \{ 0 \}$.}$$

In this paper, we will refer to a positive definite Hermitian form $h$ on $V$ more briefly as a \emph{metric} on $V$. Given a metric $h$ on $V$ the \emph{norm} $|v|$ of a vector $v \in V$ is the non-negative real number determined by
$$|v|^2=h(v,v),$$ and the \emph{adjoint} $A^h$ of $A \in \mathrm{End}(V)$ is determined by
$$h(Au,v)=h(u,A^h v) \quad \hbox{for all $u,v \in V$.}$$ The metric on $V$ then induces a metric on $\mathrm{End}(V)$ by defining 
$$h(A,B)=\mathrm{tr}(A^h B) \quad \hbox{for $A,B \in \mathrm{End}(V)$.}$$ The orthogonal direct sum of $m$ copies of $\mathrm{End}(V)$ has the metric given by
$$h(A,B)=\sum_{i=1}^m h(A_i,B_i) \quad \hbox{for $A=(A_1,\dots,A_m), B=(B_1,\dots,B_m) \in \mathrm{End}(V)^{\oplus m}$.}$$ We hope that the use of the same letter $h$ for all these will not cause confusion.

\subsection{Minimal metrics}
Let $V$ be a finite-dimensional non-zero $\CC$-vector space and let $A=(A_1,A_2,\dots,A_m)$ be a sequence of $\CC$-linear endomorphisms of $V$. We write
$$M=\{ \hbox{metrics on $V$} \} \cong K \backslash G$$ for the space of metrics on $V$, where $G=\mathrm{GL}(V)$ and $K=U(h_0)$ is the subgroup of $G$ fixing a particular choice $h_0$ of metric. The \emph{$A$-energy} of a metric $h$ is the quantity
$$E_A(h)=\sum_{i=1}^m \mathrm{tr}(A_i^h A_i).$$ Evidently $E_A(h) \geq 0$, and $E_A(h)=0$ only if all $A_i$'s are $0$. A metric $h$ is \emph{$A$-minimal} if $E_A(h)$ is minimal. The next lemma shows that if $A$ acts irreducibly on $V$, then minimal metrics exist and are unique up to multiplication by $\RR_{>0}$.

\begin{lemma} \label{minimal metric}
Suppose $0$ and $V$ are the only subspaces of $V$ invariant by all the $A_i$'s. Then the metric $h$ is minimal if and only if 
$$0=\sum_{i=1}^m [A_i^h,A_i],$$ and furthermore a minimal metric exists and is unique up to multiplication by $\RR_{>0}$. 
\end{lemma}
\begin{proof}
Fix a metric $h_0$ on $V$. The group $G$ acts by conjugation on $\mathrm{End}_\CC(V)$ and hence on $R=\mathrm{End}_\CC(V)^{\oplus m}$, as well as on the set $M$ of metrics. The sequence $A$ is an element of $R$, and by the Hilbert-Mumford criterion, the simplicity of $V$ as an $A$-module implies that the orbit $G \cdot A$ is closed in $R$. The unitary group $K=U(h_0)$ preserves the metric $h_0$ and hence the metric it induces on $R$. So by the Kempf-Ness theorem (Theorem 0.1 of \cite{KeNe}) the set of points of $G\cdot A$ of minimal distance to the origin is a single $K$-orbit, equal to the set of critical points of the function
$$E_A(g)=\sum \left| g\cdot A_i \right|^2,$$ which descends to $M$ to give the function $E_A$ on $M$. Direct computation shows that the critical points of this function $E_A$ on $G$ are the points $g \in G$ with
$$\sum_{i=1}^m [(g A_i g^{-1})^{h_0}, g A_i g^{-1}]=0,$$ which becomes the condition 
$$0=\sum_{i=1}^m [A_i^h,A_i]$$ upon conjugation by $g$, and with $h=g^{-1} \cdot h_0$. The uniqueness up to scalars follows from the fact that Kempf-Ness implies that the minima are achieved on a single $K$-orbit, together with the irreducibility of $V$ as an $A$-module.
\end{proof}

The uniqueness assertion above allows us to prove that certain transformations preserve a minimal metric:\begin{corollary} \label{inv metric}
With the hypotheses of the previous lemma, suppose $x \in G$ is such that:
\begin{itemize}
 \item[(a)] There exist $w \in S_m$ and complex numbers $\zeta_1,\dots,\zeta_m$ of absolute value $1$ such that for all $1 \leq i \leq m$,  $x A_i x^{-1}=\zeta_i A_{w(i)}$, and
 \item[(b)] some eigenvalue of $x$ has absolute value $1$.
 \end{itemize} Then if $h$ is a minimal metric for $A$, we have $h(xu,xv)=h(u,v)$, or in other words $x$ preserves the metric $h$.
\end{corollary}
\begin{proof}
This follows from the uniqueness in the previous lemma.
\end{proof}

\subsection{The rational Cherednik algebra of the symmetric group}

The rational Cherednik algebra (at $t=0$ and $c=1$) of the symmetric group is the $\CC$-algebra $H$ with generators $x_1,\dots,x_n,y_1,\dots,y_n$, and the symmetric group $S_n$, subject to the relations
$$w x_i w^{-1}=x_{w(i)} \quad w y_i w^{-1}=y_{w(i)} \quad \hbox{for all $1 \leq i \leq n$ and all $w \in S_n$,}$$
$$[x_i,x_j]=0=[y_i,y_j] \ \hbox{for all $1 \leq i,j \leq n$,} \quad [y_i,x_j]=(ij) \ \hbox{for all $1 \leq i \neq j \leq n$,} $$ and
$$[y_i,x_i]=-\sum_{j \neq i} (ij).$$

The irreducible $H$-modules are all, as $S_n$-modules, isomorphic to the regular representation (this is the $n!$ theorem for $H$, and it is easy to check thanks to an extremely elegant argument of Etingof). For each partition $\lambda$ of $n$, there is a unique irreducible $H$-module $L(\lambda)$ on which the $x$'s and $y$'s act nilpotently and for which the subspace annihilated by the $y$'s is isomorphic to $S^\lambda$ as an $S_n$-module. Filtering $H$ by total degree in the $x$'s and $y$'s, so that $F_\ell H$ is spanned by products of elements of the symmetric group by at most $\ell$ of the symbols $x_1,\dots,x_n,y_1,\dots,y_n$, we also filter $L$ by $F_\ell =(F_\ell H) \cdot L^{S_n}$. Fixing a non-zero $t \in L^{S_n}$ we also have $F_\ell=F_\ell H \cdot t$.

Thanks to Haiman's characterization of $\mathrm{GH}_\lambda$ in Propositions 4.1.2 and 4.1.3 of \cite{Hai3}, the main theorem of  \cite{Gri1} implies
\begin{theorem}
If $\mathrm{gr}(L(\lambda))$ has a one-dimensional socle, then it is isomorphic to $\mathrm{GH}_\lambda$, and so the $n!$ theorem holds for $\lambda$. 
\end{theorem} In the remainder of this article we fix $\lambda$ and write $L=L(\lambda)$. Our goal is to prove that its associated graded $\mathrm{gr}(L)$  has a one-dimensional socle. The obstacle to doing so is that it is difficult to exert control over the filtration $F_\ell$. We will do so by constructing a second grading, using a rotation of the $H$-action.

\subsection{Quaternions}

The division algebra of quaternions $\HH$ is the $\RR$-algebra with basis $1,i,j,k$ and multiplication determined by the rules
$$i^2=j^2=k^2=ijk=-1.$$ Because it is non-commutative, one must distinguish between left and right vector spaces over it. A left $\HH$-vector space $V$ is in particular a $\CC$-vector space, and multiplication by $j$ on the left induces a $\CC$-conjugate linear endomorphism of $V$ with square $-1$. Conversely, given a $\CC$-vector space $V$ equipped with a $\CC$-conjugate linear endomorphism $J$ with $J^2=-1$, we obtain a quaternionic vector space by defining multiplication by $j$ to be given by the operator $J$. 

We note that $\HH$ itself is a $\CC$-vector space with basis $1,j$, and that left multiplication by $j$ on this basis is the map
$$a+bj \mapsto -\overline{b} + \overline{a}j,$$ which may motivate the formula below for our construction of the rotated $H$-action.

We note that we will frequently use the letters $i,j$, and $k$ as indices of summation having nothing whatsoever to do with quaternions. We hope this will not cause confusion.

\subsection{Rotation of the $H$-action}

If $L=L(\lambda)$ is the irreducible $H$-module above, then the elements $x_i,y_i \in H$ determine an element $A=(x_1,\dots,x_n,y_1,\dots,y_n)$ of $E=\mathrm{End}_\CC(L)^{\oplus 2n}$, and since they generate $H$, the vector space $L$ is an irreducible $A$-module. Now $E$ is not an $\HH$-vector space in a natural fashion, but if we fix a metric $h$ on $L$ and define
$$j_h \cdot (x_1,\dots,x_n,y_1,\dots,y_n)=(-y_1^h,\dots,-y_n^h,x_1^h,\dots,x_n^h)$$ then since $j_h$ is $\CC$-conjugate linear and satisfies $j_h^2=-1$, we obtain a left $\HH$-vector space structure on $E$. This allows for the possibility of varying the operators $x_i$ and $y_i$ by multiplying them by quaternions, \emph{rotating} them. On the other hand, we can also choose a particular quaternion and then vary $h$, which is what we will do. 

A choice which seems to be maximally convenient for calculation is multiplication by $(1+j_h)/2$. Thus we define $p_i$ and $q_i$ by
$$\frac{1}{2}(1+j_h)(x_1,\dots,x_n,y_1,\dots,y_n)=(p_1,\dots,p_n,q_1,\dots,q_n),$$ or more explicitly
$$p_i=\frac{1}{2}(x_i-y_i^h) \quad \text{and} \quad q_i=\frac{1}{2}(y_i+x_i^h).$$ Of course $p_i$ and $q_i$ depend on the choice of $h$, but we leave it out of the notation, as we will soon fix a particular choice. Our first main theorem is now:

\begin{theorem}
If the metric $h$ is $A$-minimal then the operators $$p_1,\dots,p_n,q_1,\dots,q_n$$ commute pairwise. Conversely, if the metric $h$ is $S_n$-invariant and the operators commute pairwise, then $h$ is $A$-minimal.
\end{theorem}

\begin{proof}
Assume the metric $h$ is $S_n$-invariant (which holds automatically if $h$ is $A$-minimal, by Corollary \ref{inv metric}). Direct calculation (making use of the relations of $H$ and the fact that $h$ is $S_n$-invariant, so that transpositions are self-adjoint) shows that all the commutators of the $p$'s and $q$'s are zero if and only if the following relations hold:
\begin{equation}
[x_i,x_j^h]+[y_j,y_i^h]=0 \quad \text{and} \quad [x_i,y_j^h]-[x_j,y_i^h]=0 \quad \hbox{for all $1 \leq i,j \leq n$.}
\end{equation} Thus the $p$'s and $q$'s commute pairwise if and only if
$$0=\sum_{1 \leq i,j \leq n} \left| [x_i,x_j^h]+[y_j,y_i^h] \right|^2+\left|  [x_i,y_j^h]-[x_j,y_i^h] \right|^2,$$ where the norm here is that corresponding to $h$, 
$$|A|^2=\mathrm{tr}(A^h A).$$ Using the invariance $\mathrm{tr}([A,B]C)=\mathrm{tr}(A[B,C])$ and the Jacobi identity, we may rewrite
\begin{align*}
\sum_{1 \leq i,j \leq n} & \left| [x_i,x_j^h]+[y_j,y_i^h] \right|^2+\left|  [x_i,y_j^h]-[x_j,y_i^h] \right|^2 \\ 
&=\left| \sum_{i=1}^n [x_i,x_i^h]+[y_i,y_i^h] \right|^2+2 \left( \sum_{1 \leq i,j \leq n} \left| [x_i,y_j] \right|^2-\left| \sum_{i=1}^n [x_i,y_i] \right|^2 \right)
\end{align*}
Direct calculation, using that $L$ is the regular representation of $S_n$, shows that
$$\sum_{1 \leq i,j \leq n} \left| [x_i,y_j] \right|^2-\left| \sum_{i=1}^n [x_i,y_i] \right|^2=0,$$ so that commutativity of the $p$'s and $q$'s is equivalent to
$$0=\sum_{i=1}^n [x_i,x_i^h]+[y_i,y_i^h].$$ According to Lemma \ref{minimal metric}, this is equivalent to minimality of the metric $h$.
\end{proof}

\subsection{Consequences of the existence of a minimal metric}

From now on we fix an $A$-minimal metric $h=h_A$ for $A=(x_1,\dots,x_n,y_1,\dots,y_n)$ acting on $L=L(\lambda)$. First we show that the Euler grading (the unique grading for which $t$ is degree $0$, the operators $x$ have degree $1$, and the operators $y$ have degree $-1$) on $L$ is $h$-orthogonal:
\begin{lemma}
The grading subspaces $L_d$ are mutually orthogonal with respect to $h$. In particular, the $h$-adjoint of an operator of degree $k$ with respect to this grading is of degree $-k$.
\end{lemma}
\begin{proof}
We define a map $e:S^1 \to \mathrm{GL}(L)$ by setting $e(z)(v)=z^d v$ if $v \in L_d$. Then as operators on $L$ we have
$$e(z) x_i e(z)^{-1}=zx_i \quad \text{and} \quad e(z) y_i e(z)^{-1}=z^{-1} y_i.$$ Corollary \ref{inv metric} shows that this action preserves $h$, implying that the $L_d$'s are $h$-orthogonal subspaces. \end{proof}

Next we record the additional relations between $p$'s and $q$'s that we will need. The proof is a routine calculation.
\begin{lemma} \label{ps and qs}
For $1 \leq i \leq n$ and $w \in S_n$ we have
$$w p_i w^{-1}=p_{w(i)} \quad \text{and} \quad w q_i w^{-1}=q_{w(i)},$$
while for all $1 \leq i,j \leq n$ we have
$$[p_i,p_j^h]+[q_j,q_i^h]=[y_i,x_j] \quad \text{and} \quad  [p_i,q_j^h]=[p_j,q_i^h].$$ \end{lemma}

We have $x_i=p_i+q_i^h$ and $y_i=q_i-p_i^h$. Motivated by these formulas, for $z \in S^1$ we define the elements
$$x_i^{(z)}=zp_i+z^{-1} q_i^h \quad \text{and} \quad y_i^{(z)}=zq_i-z^{-1} p_i^h.$$ These provide an $S^1$-worth of $H$-actions on $L$, which together with the rigidity of $L$ gives a second grading.

\begin{lemma}
These define a representation of $H$ on $L$, which we denote $L^{(z)}$. Moreover $L^{(z)} \cong L$ as $H$-modules.
\end{lemma}
\begin{proof}
Direct calculation shows that this construction produces a representation $L^{(z)}$ of $H$. Since it is $n!$ dimensional it is automatically irreducible, the new $x$'s and $y$'s have Euler degree $1$ and $-1$ and so act nilpotently on it, and the space of singular vectors it contains is $S^\lambda$ so it is isomorphic to $L(\lambda)$.
\end{proof}
\subsection{The second grading on $L$} 

\begin{lemma}
We let $U(h) \leq G$ be the unitary group of the form $h$, consider the conditions 
\begin{equation} \label{g rels} g x_i g^{-1}=x_i^{(z)}, \quad g y_i g^{-1}=y_i^{(z)} \quad \text{and} \quad g w g^{-1}=w \quad \hbox{for all $1 \leq i \leq n$ and all $w \in S_n$.} \end{equation}
and define 
$$\Gamma=\{(z,g) \ | \ z \in S^1, \ g \in U(h), \ gt=t, \ \hbox{and $g$ satisfies \eqref{g rels}} \}.$$ The projection of $\Gamma$ onto its first factor is a Lie group isomorphism to $S^1$, and composing its inverse with the projection of $\Gamma$ onto its second factor gives a Lie group homomorphism from $S^1$ to $U(h)$.
\end{lemma}

\begin{proof}
First observe that the conditions defining $\Gamma$ are closed conditions, so it is a compact subset of $S^1 \times U(h)$. Moreover it is a subgroup: take $h$-adjoints in the equation $g y_i g^{-1}=y_i^{(z)}$ and use that $g$ is unitary to conclude that $gp_i g^{-1}=zp_i$ and $gq_i g^{-1}=zq_i$. It thus suffices to prove that the projection to $S^1$ is bijective, for it is then a homeomorphism.

The projection of $\Gamma$ on $S^1$ is injective, since its kernel consists of pairs $(1,g)$ where $g$ is an $H$-module endomorphism of $L$ fixing $t$, and hence equal to the identity.

We show the projection is surjective. The previous lemma implies that for each $z \in S^1$, there is some $g \in \mathrm{GL}(L)$ satisfying \eqref{g rels}. Since $t=t_\lambda$ is a basis of $L^{S_n}$, after multiplying $g$ by a scalar we may assume $g(t)=t$. It remains to show that $g \in U(h)$. We compute, using that $p$'s and $q$'s commute with one another together with Lemma \ref{ps and qs}
\begin{align*}
\sum_{i=1}^n & [x_i^{(z)},(x_i^{(z)})^h]+[y_i^{(z)},(y_i^{(z)})^h]=\sum_{i=1}^n [z p_i+z^{-1}q_i^h,z^{-1} p_i^h+z q_i]+[zq_i-z^{-1} p_i^h,z^{-1} q_i^h-z p_i] \\
&=\sum_{i=1}^n [p_i,p_i^h]+[q_i^h,q_i]+[q_i,q_i^h]+[p_i^h,p_i]=\sum_{i=1}^n [y_i,x_i]-[y_i,x_i]=0.
\end{align*} It thus follows from the uniqueness of minimal metrics that there is $c \in \RR_{>0}$ such that $h=c g \cdot h$. But $0<h(t,t)=c h(g^{-1} t,g^{-1}t)=c h(t,t)$ then gives $c=1$, so $h$ is $g$-invariant. \end{proof} We write $z \mapsto g(z)$ for the homomorphism provided by the lemma. 

Let $L^{d}=\{v \in L \ | \ g(z)v=z^d v \ \hbox{for all $z \in S^1$}\}.$ Then $L$ is the direct sum $L=\bigoplus_{d \in \ZZ} L^d$ and the operators $p$ and $q$ are of degree $1$ with respect to this grading by the first paragraph of the proof of the lemma. 

\section{The structure of $L$ as a $\CC[p,q]$-module and the proof of the $n!$ theorem}

\subsection{The structure of $L$ as a module over $\CC[p_1,\dots,p_n,q_1,\dots,q_n]$}

We abbreviate by writing $\CC[p,q]=\CC[p_1,\dots,p_n,q_1,\dots,q_n]$.

\begin{lemma}
We have
$$L^{S_n}=\{v \in L \ | \ p_i^hv=0=q_i^hjv \ \hbox{for all $1 \leq i \leq n$} \} \quad \text{and} \quad L^\mathrm{det}=\{v \in L \ | \ p_i v=0=q_i v \ \hbox{for all $1 \leq i \leq n$.} \}.$$ Moreover, $t$ generates $L$ as a $\CC[p,q]$-module, and the socle of this $\CC[p,q]$-module is $L^\mathrm{det}$.
\end{lemma} 
\begin{proof}
The operators $p_i$ and $q_i$ are nilpotent and commute with one another, so the space of vectors they annihilate is non-zero. Likewise for $p_i^h$ and $q_i^h$. Fix a subset $S \subseteq \{1,2,\dots,n \}$ and put $p_S=\sum_{i \in S} p_i$ and $q_S=\sum_{i \in S} q_i$. Then $$[p_S,p_S^h]+[q_S,q_S^h]=\sum_{i,j \in S}  ([p_i,p_j^h]+[q_i,q_j^h])=\sum_{ i,j \in S} [y_i,x_j]=-\sum_{\substack{i \in S \\ j \notin S}} (ij).$$ Now if $p_i^hu=0=q_i^hu$ for all $1 \leq i \leq n$ then we obtain
$$h(([p_S,p_S^h]+[q_S,q_S^h])u,u)=-h((p_S^h p_S+q_S^h q_S )u,u)=-h(p_S u,p_S u)-h(q_S u,q_S u) \leq 0.$$ A calculation with contents shows that for any $\mu \neq (n)$, taking $S=\{1,2,\dots,m\}$ where $m$ is the number of boxes in $\mu$ without its last row, the element $-\sum_{\substack{i \in S \\ j \notin S}} (ij)$ has a positive eigenvalue on $S^\mu$. Hence $L^{S_n}$ is the space annihilated by all the $p_i^h$'s and $q_i^h$'s. Similarly, $L^{\mathrm{det}}$ is the subspace annihilated by all $p_i$'s and $q_i$'s. The first two assertions of the lemma are proved, as is the statement about the socle.

It follows that $(L^{S_n})^\perp \supseteq \sum_{i=1}^n p_i L + q_i L$, and likewise $$\left(\sum_{i=1}^n p_i L + q_i L \right)^\perp \subseteq  \{v \in L \ | \ p_i^*v=0=q_i^*v \ \hbox{for all $1 \leq i \leq n$} \}=L^{S_n}.$$ Hence $L=L^{S_n} \oplus \sum_{i=1}^n p_i L + q_i L$, from which the remaining assertion of the lemma follows.
\end{proof}

The lemma implies that $L$ is in fact bigraded by $L=\bigoplus L_d^e$, where $L_d^e=L_d \cap L^e$. 
\subsection{The filtration $F_\ell$} 

Now we obtain the promised control over $F_\ell$ and prove the $n!$ theorem. 
\begin{theorem}
For each $\ell \in \ZZ_{\geq 0}$ we have $F_\ell=\bigoplus_{d=0}^\ell L^d$, and there is an isomorphism $\mathrm{gr}(L) \cong L$ of $\CC[p,q]$-modules, where $p_i$ and $q_i$ act on $\mathrm{gr}(L)$ as $x_i$ and $y_i$.
\end{theorem} 
\begin{proof}
We recall that the grading $L^d$ is defined by $$L^d=\{v \in L \ | \ g(z)v=z^d v \ \hbox{for all $z \in S^1$}\}$$ and that $p_i$ and $q_i$ are degree $1$ for this grading. Moreover, $t$ is fixed by $g(z)$ implies that $L^{S_n}$ is in degree $0$, and is equal to $L^0$ since $L=\CC[p,q] t$. This last equation shows more generally that $L^{\leq \ell}=\bigoplus_{d=0}^\ell L^d$ is the span of all polynomials of degree at most $\ell$ in $p$'s and $q$'s applied to $t$. But $x_i=p_i+q_i^h$ and $y_i=q_i-p_i^h$ map $L^{\leq \ell}$ into $L^{\leq \ell+1}$ so we obtain the containment $L^{\leq \ell} \supseteq F_\ell$. The other containment is also established by induction on $\ell$ using that for any polynomial $f$ of degree $\ell$, and fixing an ordering of the variables $x$ and $y$, 
$$f(x,y)t -f(p,q)t=f(p+q^h,q-p^h)t-f(p,q)t \in L^{\leq \ell-1}.$$ For the same reason, the last assertion of the theorem follows.
\end{proof} This theorem together with the preceding lemma proves that $\mathrm{gr}(L)$ has a one-dimensional socle, from which the $n!$ theorem follows.

\def\cprime{$'$} \def\cprime{$'$}

\end{document}